\pdfoutput=1
\RequirePackage{ifpdf}
\ifpdf 
\documentclass[pdftex]{sigma}
\else
\documentclass{sigma}
\fi

\newtheorem{thm}{Theorem}[section]
\newtheorem{lem}[thm]{Lemma}

\newtheorem{prop}[thm]{Proposition}
{ \theoremstyle{definition}
\newtheorem{Example}[thm]{Example}
\newtheorem{Remark}[thm]{Remark} }

\numberwithin{equation}{section}

\begin{document}


\renewcommand{\PaperNumber}{061}

\FirstPageHeading

\ShortArticleName{Spectral Analysis of Certain Schr\"{o}dinger Operators}

\ArticleName{Spectral Analysis of Certain Schr\"{o}dinger Operators}

\Author{Mourad E.H.~ISMAIL~$^\dag$ and Erik KOELINK~$^\ddag$}

\AuthorNameForHeading{M.E.H.~Ismail and E.~Koelink}

\Address{$^\dag$~Department of Mathematics, University of Central Florida, Orlando, FL 32816, USA}
\EmailD{\href{mailto:mourad.ismail@ucf.edu}{mourad.ismail@ucf.edu}}
\URLaddressD{\url{http://www.math.ucf.edu/~ismail/}}

\Address{$^\ddag$~Radboud Universiteit, IMAPP, FNWI, Heyendaalseweg 135, 6525 AJ Nijmegen, \\
\hphantom{$^\ddag$}~the Netherlands}
\EmailD{\href{mailto:e.koelink@math.ru.nl}{e.koelink@math.ru.nl}}
\URLaddressD{\url{http://www.math.ru.nl/~koelink/}}

\ArticleDates{Received May 07, 2012, in f\/inal form September 12, 2012; Published online September 15, 2012}

\Abstract{The $J$-matrix method is extended to dif\/ference and $q$-dif\/ference operators and is applied to
several  explicit dif\/ferential, dif\/ference, $q$-dif\/ference and second order Askey--Wilson type operators. The
spectrum and the spectral measures are discussed in each case and the corresponding eigenfunction
expansion is written down explicitly
in most cases. In some cases we encounter new orthogonal polynomials with explicit three term
recurrence relations where nothing is known about their explicit representations or
orthogonality measures. Each model we analyze is a discrete quantum mechanical model in the
sense of Odake and Sasaki~[\textit{J.~Phys.~A: Math.
  Theor.} \textbf{44} (2011), 353001, 47~pages].}


\Classification{30E05; 33C45; 39A10; 42C05; 44A60}

\section{Introduction}\label{section1}

The $J$-matrix method started with the pioneering works of Yamani,  Heller
and Reinhardt \cite{Hel,Hel:Rei:Yam, Yam:Rei} in the early 1970's and has been
applied by Yamani,  Heller and Reinhardt to dif\/ferent physical models. Some of the recent
applications of the $J$-matrix method to physics are spearheaded by Alhaidari and his research
team,~\cite{AlH4, AlH2,AlH1,AlH3,AlHaBAAA}.  The $J$-matrix principle
 says that the spectrum of a~tridigonalizable operator is the same as the tridiagonal matrix
 representing it. Such a~tridiagonal matrix can be split into irreducible blocks, and to each of
 these blocks there is a~corresponding set of orthogonal polynomials. Moreover, the eigenfunctions
 of a~tridiagonalizable operator can be expressed using these corresponding polynomials, and in the
 self-adjoint case the spectral measure is related to the orthogonality measures of the orthogonal
 polynomials. This general set-up is described and proved in~\cite{Ism:Koe} where we considered the
 Schr\"odinger equation with the Morse potential as an example.  Our later work~\cite{Ism:Koe2}
 develops a general scheme for tridiagonalizing dif\/ferential, dif\/ference, or $q$-dif\/ference
 operators arising from two sets of related orthogonal polynomials. In particular, we f\/ind in~\cite{Ism:Koe2} the
 Jacobi transform and its special case the Mehler--Fock transform  originally introduced by Mehler
 to study electrical distributions.

Tridiagonalization of an explicit symmetric or self-adjoint
operator $T$, like a dif\/ferential or \mbox{($q$-)}dif\/ference operator
on an explicit Hilbert space $\mathcal{H}$ of functions, amounts to f\/inding an explicit orthonormal
basis $\{e_n\}_{n=0}^\infty$ so that we can realize
\[
Te_n  =   \alpha_n e_{n+1} + \beta_ne_n + \alpha_{n-1}e_n, \qquad \alpha_n>0, \quad \beta_n \in {\mathbb{R}},
\]
with respect to this basis. Since three-term operators, or Jacobi operators, can be described
in terms of orthogonal polynomials, we can obtain information on $T$ in terms of properties of
the orthogonal polynomials, see \cite{Ism:Koe} and references given there. Assuming that there
are orthgonal polynomials $p_n$ (corresponding to a determinate moment problem) satisfying
\[
\lambda p_n(\lambda)   =   \alpha_n p_{n+1}(\lambda) + \beta_np_n(\lambda) + \alpha_{n-1}p_n(\lambda),
\qquad \int_{\mathbb{R}} p_n(\lambda) p_m(\lambda) \, d\mu(\lambda) = \delta_{nm},
\]
then the spectral decomposition is $UT=MU$, where $U\colon \mathcal{H}\to L_2(\mu)$, $e_n\mapsto p_n$ and
$M\colon L_2(\mu)\to L_2(\mu)$ is the multiplication operator.
The link between more general dif\/ferential, dif\/ference or $q$-dif\/ference operators, and Jacobi
matrices can be very useful to study the spectral decomposition of the original operator in terms
of the orthogonal polynomials and vice versa. So we can obtain information on one of the operators
by transferring information from the other, and we show this in particular examples in this paper.
In particular, in case information on both operators is known, we obtain even more explicit
results, and examples of this approach can be found in~\cite{Ism:Koe2}. However, it is not
straightforward to f\/ind explicit tridiagonalization of a given operator, and here we present
ways to obtain operators with a tridiagonalization.
In the tridiagonalization
of the operators in this paper it is often the case that the polynomials cannot be matched directly
with known polynomials~\cite{Ism, Koe:Swa}, and in these cases we have given some information
on the support of the spectral meaure of these polynomials.
A treatment of  the spectral theory of dif\/ferential operators can be found in many sources, and we
refer the interested reader to the excellent book by Edmunds and Evans~\cite{Edm:Eva}. The spectral
theory of tridiagonal matrices and their connection with orthogonal polynomials and the moment
problem is in~\cite{Akh, Deift, Koel, Sim}.

\looseness=-1
The contents of the paper are as follows.
In Section~\ref{section2} we record the formulas used in the sequel. The basic references are
\cite{And:Ask:Roy,Erd:Mag:Obe:Tri1,Gas:Rah,Ism,Koe:Swa,Rai}.  The expert  reader may easily skip
Section~\ref{section2}.  In Sections~\ref{section3} and~\ref{section4} we start with an operator with known orthogonal polynomial
eigenfunctions then multiply it by a linear function of the space variable and study the spectral
properties of the new Hamiltonian. Section~\ref{section3} treats the case of Laguerre polynomials, leading to
tridiagona\-li\-zation involving continuous dual Hahn polynomials. It is simple enough but contains all
the ingredients of the method. Section~\ref{section4} treats the Meixner polynomials and the $J$-matrix method
leads to a one parameter generalization of the continuous dual Hahn polynomials. The examples in
Sections~\ref{section3} and~\ref{section4} are related to the approach of~\cite{Ism:Koe2}.
In Section~\ref{section5} we introduce a~dif\/ferent modif\/ication. We start with an operator $T$ which is
diagonalized by a known polynomial system.  We then consider the Sch\"{o}dinger operator
$T+\gamma x$.  Such an equation arises for example in the case of a charged particle in the
presence of a uniform electric f\/ield. In this case $\gamma = -qF$, where $F$ is the magnitude of
the electric f\/ield and $q$ is the electric charge,
see \cite[\S~24]{Lan:Lif}, or~\cite{Bas}.

The Sch\"{o}dinger operator $T+\gamma x$  is automatically tridiagonal in the orthogonal polynomial
basis which
diagonalizes or tridiagonalizes
$T$. We study the spectral properties of $T+\gamma x$ for four dif\/ferent sets of
polynomials. In Section~\ref{section6} we combine both generalizations of Sections~\ref{section3} and~\ref{section4} for the
case of the Laguerre and Meixner polynomials.
Finally, in Section~\ref{section7} we study this approach for two families of $q$-orthogonal polynomials, namely
the Al-Salam--Chihara polynomials and the $q^{-1}$-Hermite polynomials. The $q^{-1}$-Hermite
polynomials correspond to an indeterminate moment problem, so we study the corresponding $q$-dif\/ference operator on the weighted $L^2$-space corresponding to a $N$-extremal measure. It turns out
that the polynomials in the tridiagonali\-za\-tion is again corresponding to an indeterminate moment
problem, so that the $q$-dif\/ference operator is not essentially self-adjoint on the space of
polynomials.

We end by noting that more dif\/ferential, dif\/ference and $q$-dif\/ference operators
can be studied using the $J$-matrix method. In particular, we can study classes of higher-order
operators in this way as well.

\section{Preliminaries}\label{section2}

In this section we recall some results needed in the sequel.
 We f\/irst record the properties of the Laguerre polynomials.   They   satisfy the dif\/ferential
 relations,
 \cite[(4.6.13), (4.6.15)]{Ism}:
 \begin{gather}
\frac{d}{dx}   L_n^{(\alpha)}(x) = -L_{n-1}^{(\alpha+1)}(x), \label{eqLaglower}
\\
\frac{1}{x^\alpha e^{-x}}  \frac{d}{dx}  \left[x^{\alpha+1} e^{-x} \frac{d}{dx}  L_n^{(\alpha)}
 (x)\right]
  = -n L_n^{(\alpha)}(x).  \label{eqLag2ndorder}
 \end{gather}
A generating function  of the Laguerre polynomials is
 \begin{gather*}
 \sum_{n=0}^\infty L_n^{(\alpha)}(x) t^n = (1-t)^{-\alpha-1}\exp\left(\frac{-xt}{1-t}\right),
 \end{gather*}
see \cite[(4.6.4)]{Ism}, \cite{Koe:Swa, Rai} and it implies
 \begin{gather}
 \label{eqLaguerrealandalplus1}
 L_n^{(\alpha)}(x) = L_n^{(\alpha+1)}(x) - L_{n-1}^{(\alpha+1)}(x).
  \end{gather}
  The orthogonality relation is
  \begin{gather}
  \int_0^\infty x^\alpha e^{-x}  L_m^{(\alpha)}(x) L_n^{(\alpha)}(x) dx = \frac{\Gamma(\alpha+n+1)}{n!} \, \delta_{m,n}, \qquad \alpha >-1.
  \label{eqLagorth}
  \end{gather}

  The Meixner polynomials are,  \cite[\S~6.1]{Ism}, \cite[\S~1.9]{Koe:Swa},
  \begin{gather*}
   M_n(x;\beta,c) =  {}_2F_1(-n, -x; \beta; 1-1/c),
  \end{gather*}
  and have the generating function
  \begin{gather}
  \label{eqGFMeixner}
  \sum_{n=0}^\infty \frac{(\beta)_n}{n!} M_n(x;\beta,c) t^n = (1-t/c)^x(1-t)^{-x-\beta}.
  \end{gather}
  The orthogonality relation is
  \begin{gather*}
  \sum_{x=0}^\infty  M_m(x;\beta,c)  M_n(x;\beta,c) \frac{(\beta)_x}{x!} c^x = \frac{c^{-n}\, n!}{(\beta)_n(1-c)^\beta} \delta_{m,n},
  \end{gather*}
  valid for $\beta >0$, $0<c<1$ and their three term recurrence relation is
 \begin{gather}
 -x M_n(x;\beta,c) =\frac{c(\beta+n)}{1-c}M_{n+1}(x;\beta,c) + \frac{n}{1-c}M_{n-1}(x;\beta,c)
\nonumber\\
\hphantom{-x M_n(x;\beta,c) =}{}
- \frac{n+c(\beta+n)}{1-c}M_n(x;\beta,c). \label{eqMei3trr}
 \end{gather}
 The Meixner polynomials satisfy the second order dif\/ference equation
  \begin{gather}
  \label{eqDEMeixner}
 \frac{x!}{c^x(\beta)_x} \nabla\left(\frac{(\beta+1)_xc^x}{x!}\Delta\right)M_n(x;\beta,c) = \frac{n}{\beta}
 \frac{c-1}{c} M_n(x;\beta,c),
  \end{gather}
where
\begin{gather*}
(\Delta f)(x) = f(x+1)-f(x), \qquad (\nabla f)(x) = f(x)-f(x-1).
\end{gather*}

The Meixner--Pollaczek polynomials $\{P_n^{(\lambda)}\!(x;\phi)\}$ satisfy the orthogonality relation
 \cite[(1.7.2)]{Koe:Swa}
\begin{gather*}
\frac{1}{2\pi} \int_{\mathbb{R}}  e^{(2\phi-\pi)x}|\Gamma(\lambda+ix)|^2 P_m^{(\lambda)}(x;\phi)P_n^{(\lambda)}(x;\phi) \, dx = \frac{\Gamma(n+2\lambda)}{(2\sin \phi)^{2\lambda}  n!} \delta_{m,n},
\end{gather*}
for $\lambda >0$, $0< \phi < \pi$, and the three term recurrence relation  \cite[(1.7.3)]{Koe:Swa}
\begin{gather}
\label{eqMeixPoll3trr}
(n+1)P_{n+1}^{(\lambda)}(x;\phi) + (n+2\lambda -1)P_{n-1}^{(\lambda)}(x;\phi)
= 2[x\sin \phi +(n+\lambda)\cos \phi]P_n^{(\lambda)}(x;\phi),
\end{gather}
with $P_0^{(\lambda)}(x;\phi) =1$, $P_1^{(\lambda)}(x;\phi)= 2[x\sin \phi + \lambda \cos \phi]$.

 We parametrize  the independent variable $x$ by $ x = (z+1/z)/2$ and given a function  we set
 $\breve{f}(z) = f(x)$.
 The Askey--Wilson operator $\mathcal{D}_q$ and the averaging operator  $\mathcal{A}_q$ are def\/ined by, \cite[\S~12.1]{Ism},
 \begin{gather*}
 (\mathcal{D}_q f)(x) = \frac{\breve{f}(zq^{1/2})-  \breve{f}(zq^{-1/2})}{ \breve{e}(q^{1/2}z)- \breve{e}(q^{-1/2}z)}, \qquad
(\mathcal{A}_q f)(x) = \frac{1}{2}\big[ \breve{f}(zq^{1/2})+  \breve{f}(zq^{-1/2})\big],
 \end{gather*}
 where $e(x) = x = (z+1/z)/2$.

The Askey--Wilson operator  is well-def\/ined on $H_{1/2}$, where
\begin{gather*}
H_{\nu}:=\left\{f: f((z+1/z)/2)\text{ is analytic for } q^{\nu}\le |z| \le q^{-\nu}\right\}.
\end{gather*}
Let ${\cal H}_w$ denote the weighted space $L^2(-1,1;w(x)dx)$ with inner product
\begin{gather*}
(f,g)_w:=\int_{-1}^1 f(x)  \overline{g(x)}  w(x)\,dx,
\qquad \|f\|_w:=(f,f)_w^{1/2}
\end{gather*}
and let $T$ be def\/ined by
\begin{gather*}
Tf(x):=
 -\frac{1}{w(x)} {\cal D}_q \left(p{\cal D}_q f\right)(x),
\end{gather*}
for $f$ in $H_1$.
We shall assume that $p$ and $w$ are positive on $(-1,1)$ and also satisfy
\begin{gather*}
 \textup{(i)} \ \  p(x)/\sqrt{1-x^2} \in H_{1/2}  , \  1/p \in L(-1,1),  \nonumber\\
 \textup{(ii)} \   w(x) \in L(-1,1), \ 1/w \in L\left(-1,1; \frac{dx}{\left(1-x^2\right)}\right).
\end{gather*}
The expression $Tf$ is therefore def\/ined for $f\in H_1$, and the operator $T$ acts in $\mathcal{H}_w$.
Furthermore, the domain $H_1$ of $T$ is dense in $\mathcal{H}_w$ since it contains all polynomials. The following theorem is due to Brown, Evans and Ismail \cite{Bro:Eva:Ism}.

\begin{thm}\label{thm16.2.2}
The operator $T$ is symmetric in ${\cal H}_w$ and positive.
\end{thm}

The Al-Salam--Chihara polynomials are def\/ined by \cite[(15.1.6)]{Ism}
\begin{gather*}
p_n\left(x; t_1, t_2\,|\,q\right)= {}_{3}\phi_2\left(\left.\begin{matrix}
q^{-n}, t_1e^{i\theta}, t_1e^{-i\theta} \\
t_1t_2, 0
\end{matrix}\, \right|q,q\right).
\end{gather*}
Their weight function is
\begin{gather}
\label{eqweightAC}
w(\cos \theta ;t_1,t_2) := \frac{(e^{2i\theta}, e^{-2i\theta};q)_\infty /\sin \theta}
 { (t_1e^{i\theta}, t_1e^{-i\theta}, t_2e^{i\theta}, t_2e^{-i\theta};q)_\infty},
 \end{gather}
and their orthogonality relation will be stated in~\eqref{eqorthAC}.  The generating function for the Al-Salam--Chihara polynomials is \cite[(15.1.10)]{Ism}
\begin{gather}
\label{eqgfAC}
\sum_{n=0}^\infty \frac{(t_1t_2;q)_n}{(q;q)_n} \left(\frac{t}{t_1}\right)^n p_n(\cos \theta; t_1,t_2) =
\frac{(tt_1,tt_2;q)_\infty}{(t_1e^{i\theta}, t_1e^{-i\theta};q)_\infty}.
\end{gather}

\begin{thm}
Consider   the  three term recurrence relation in orthonormal form
\begin{gather}
\label{eqorthnormal}
xp_n(x) = a_{n+1} p_{n+1}(x) + b_np_n(x) + a_{n}p_{n-1}(x), \qquad n \ge 0,  \quad a_n > 0,
\quad b_n \in {\mathbb{R}},
\end{gather}
with $a_0 p_{-1}(x) :=0$.
Then the moment problem is determinate,  that is, it has a unique solution, if  one of the following conditions  hold
\begin{gather}
  \sum_{n=0}^\infty  \frac{|b_{n+1}|}{a_{n+1} a_{n+2}} =  \infty, \label{eqi}\\
  a_{n}+ b_n + a_{n+1}  \le  C, \qquad \text{for some} \  C,  \label{eqii}\\
  a_{n}- b_n + a_{n+1}  \le  C, \qquad \text{for some} \  C. \label{eqiii}
\end{gather}
\end{thm}

The condition \eqref{eqi} is Exercise~2 on p.~25 of \cite{Akh}, while \eqref{eqii}, \eqref{eqiii}
are Theorem~VII.1.4  and its corollary in~\cite[pp. 505--506]{Ber}.

\begin{thm}\label{IsmLi}
Let $p_n(x)$ be generated by \eqref{eqorthnormal}. Then the zeros of the polynomial $p_n(x)$  are in
$(A,B)$, where
\begin{gather*}
B = \max\{x_j: 0 < j < n\}, \qquad  A =  \min \{y_j: 0 < j < n\},
\end{gather*}
where $y_j \le x_j$ and
\begin{gather}
\label{eqdefxiyi}
x_j, y_j = \frac{1}{2} (b_j+ b_{j-1})  \pm \frac{1}{2}   \sqrt{(b_j-b_{j-1})^2 + 16 a_j^2}, \qquad 1 \le j < n.
\end{gather}
\end{thm}

Theorem \ref{IsmLi} is the  special case $c_n =1/4$  of a result due to Ismail and Li in \cite{Ism:Li}. The full result is also stated and proved in  \cite[Theorem~7.2.7]{Ism}.

The zeros of orthogonal polynomials are real and simple, so we shall follow the standard notation in~\cite{Sze} or~\cite{Ism} and arrange the zeros $x_{n,j}$, $1 \le j \le n$ as
\begin{gather*}
x_{n,1} > x_{n,2} > \cdots > x_{n,n}.
\end{gather*}

\section{A dif\/ferential operator related to the Laguerre polynomials}\label{section3}
Consider the dif\/ferential operator
\begin{gather}
\label{eq3.1}
(T_Lf) (x) = \frac{1}{x^\alpha e^{-x}}  \frac{d}{dx}  \left[x^{\alpha+2} e^{-x} \frac{df}{dx}\right].
\end{gather}
We will discuss a generalization of this operator in Section~\ref{section6}. The boundary value problem we are interested in is $T_Ly = \lambda y$ with the  boundary conditions $x^{1+\alpha/2}f(x) e^{-x/2} \to 0$ as $x \to 0$ and $x \to \infty$.  The equation $T_Ly = \lambda y$ is
\begin{gather*}
 x^2y^{\prime\prime} + (\alpha+2) x y^{\prime} - x^2 y^{\prime} = \lambda y.
\end{gather*}
It is easy to see that $T_L$ is symmetric on weighted $L_2$ space with the inner product
\begin{gather}
(f,g) = \int_0^\infty x^\alpha e^{-x} f(x) \overline{g(x)} \, dx.
\label{eqinnerprodL}
\end{gather}

The $(m,n)$ matrix elements of $T_L$ as an operator in $L_2(0, \infty, x^\alpha e^{-x})$  in the basis   $\{L_n^{(\alpha)}(x)\}$ can be calculated using \eqref{eqLaguerrealandalplus1}, \eqref{eqLagorth},  \eqref{eqLag2ndorder};
\begin{gather*}
(T_L L_n^{(\alpha)}, L_m^{(\alpha)}) = \int_0^\infty L_m^{(\alpha)}(x) \frac{d}{dx}  \left[x^{\alpha+2} e^{-x} \frac{d}{dx} L_n^{(\alpha)}(x)\right] dx  \\
=   \int_0^\infty [L_m^{(\alpha+1)}(x) -L_{m-1}^{(\alpha+1)}(x)] \frac{d}{dx} \left[x^{\alpha+2} e^{-x} \frac{d}{dx}(L_n^{(\alpha+1)}(x) -L_{n-1}^{(\alpha+1)}(x))\right] \, dx \\
=  - \int_0^\infty x^{\alpha+1} e^{-x} [L_m^{(\alpha+1)}(x) -L_{m-1}^{(\alpha+1)}(x)]
[nL_{n}^{(\alpha+1)}(x) -( n-1)L_{n-1}^{(\alpha+1)}(x)]dx  \\
=  - \frac{\Gamma(n+\alpha+2)}{(n-1)!} \delta_{m,n}  +  \frac{\Gamma(m+\alpha+2)}{(m-1)!} \delta_{m+1,n}     + \frac{\Gamma(n+\alpha+2)}{(n-1)!} \delta_{m,n+1} - \frac{\Gamma(n+\alpha+1)}{(n-2)!}\delta_{m,n} \\
= \frac{\Gamma(m+\alpha+2)}{(m-1)!} \delta_{m+1,n} -  (\alpha+2m) \frac{\Gamma(n+\alpha+1)}{(m-1)!} \delta_{m,n} + \frac{\Gamma(m+\alpha+1)}{(m-2)!} \delta_{m,n+1}.
\end{gather*}
Thus the  sought matrix representation of $T_L$ is tridiagonal. It is also clear the constants are
in the null space of $T_L$, so we mod out by  the constant functions. Let $\{A_{m,n}(L)\}$ be the
matrix elements. Thus  $A_{m,n}(L)$ is $\sqrt{\frac{m!n!}{\Gamma(m+\alpha+1)\Gamma(n+\alpha+1)}}$ times
the above expression. Thus
\begin{gather}
A_{m,n}(L) =  m\sqrt{(m+1)(m+\alpha+1)}   \delta_{m+1,n}\nonumber\\
\hphantom{A_{m,n}(L) =}{} - m(2m+\alpha)     \delta_{m,n} +
(m-1)\sqrt{m(\alpha+m)}  \delta_{m,n+1}.
\label{eq*}
\end{gather}
The ef\/fect of modding out by the constants is to delete the f\/irst row and column of the matrix is to shift $m$ and $n$ by one.  Thus we consider the tridiagonal matrix $B= (B_{m,n})$, $m,n =0, 1, \dots$,
\begin{gather*}
B_{m,n} =  (m+1)\sqrt{(m+2)(m+\alpha+2)}   \delta_{m,n-1 }  \\
\phantom{B_{m,n} =}{} -  (m+1)(2m+\alpha+2)      \delta_{m,n} +
m\sqrt{(m+1)(\alpha+m+1)}   \delta_{m,n+1}.
\end{gather*}
Now the spectral equation $E X = B X$ where $X$ is a column vector, when written componentwise  is a three term recurrence relations and the component of $X$ are
$p_n(E)$. The corresponding monic polynomials satisfy the three term recurrence relation
\begin{gather*}
Ep_m(E) = p_{m+1}(E) - (m+1)(2m+\alpha+2) p_m(E)
+ m^2(m+1)(m+\alpha+1)p_{m-1}(E).
\end{gather*}
This is \cite[(1.3.5)]{Koe:Swa} and identif\/ies the $p_m$'s as continuous dual Hahn polynomials with the parameter and variable identif\/ications
\begin{gather*}
a = \frac{1-\alpha}{2}, \qquad b = \frac{1+\alpha}{2}, \qquad c = \frac{3+\alpha}{2}, \qquad E = -x - \frac{(\alpha+1)^2}{4}.
\end{gather*}

As stated in the Introduction the spectral decomposition is given in terms of the orthogonality measure
for the continuous dual Hahn polynomials $S_n(x;a,b,c)$, see \cite[\S~1.3]{Koe:Swa}. The
measure is absolutely continuous on $[0,\infty)$ in case $a$, $b$ and $c$ are positive. In case one of
them is negative, f\/initely many discrete mass points have to be added.
In the present case however we assumed $\alpha > -1$, hence only $a$ can be negative.
Explicitly, with $S_n(x)=S_n(x;a,b,c)$, {\samepage
\begin{gather}
 \frac{1}{2\pi \Gamma(a+b) \Gamma(a+c) \Gamma(b+c)}
\int_0^\infty S_n(y^2) S_m(y^2)
\left|\frac{\Gamma(a+iy) \Gamma(b+iy) \Gamma(c+iy)}{\Gamma(2iy)}\right|^2\, dy  \nonumber \\
{} + \frac{\Gamma(b-a)\Gamma(c-a)}{\Gamma(-2a)\Gamma(b+c)} \sum_{k=0}^M
S_n\big({-}(a+k)^2\big) S_m\big({-}(a+k)^2\big)
 \frac{ (2a)_k(a+1)_k(a+b)_k(a+c)_k }{k! (a)_k(a-b+1)_k(a-c+1)_k} (-1)^k \nonumber\\  =
\delta_{nm} n!  (a+b)_n(a+c)_n(b+c)_n,\label{eq:orthorelcontHahn}
\end{gather}
where $M=\max\{k\in{\mathbb{N}}: k +a < 0\}$.}

Now \cite[Theorem~2.7]{Ism:Koe} gives the following
proposition.

\begin{prop} \label{prop:Laguerecase}\sloppy
The unbounded operator $T_L$ acting on the subspace of polynomials in
  $L_2(0, \infty, x^\alpha e^{-x})$ is essentially self-adjoint. The spectrum of the closure has
an absolutely conti\-nuous part $(-\infty, -\frac14(\alpha+1)^2]$. The discrete spectrum consists
of $\{0\}$ and $\{E_k\mid k\in \{0,\dots, M\}\}$, where $M = \max \{k\in{\mathbb{N}}\mid  k +(1-\alpha)/2 < 0\}$
and $E_k = (k+1)(k-\alpha)$.
\end{prop}

The explicit spectral measure can be obtained from the orthogonality measure for the
continuous dual Hahn polynomials~\eqref{eq:orthorelcontHahn}, cf.\ Section~\ref{section1}.

This discussion of the dif\/ferential operator $T_L$ is related to the set-up of \cite{Ism:Koe2},
where
the case related to Jacobi polynomials is considered. In \cite{Ism:Koe2} we assume that we did not
need to mod out a null space. The dif\/ferential operator $T_L$ can be related to the conf\/luent
hypergeometric dif\/ferential equation in the same way as the hypergeometric dif\/ferential operator
shows up in  \cite[\S~3]{Ism:Koe2}.

\section{A dif\/ference operator related to Meixner polynomials}\label{section4}

The generating function \eqref{eqGFMeixner} implies
\begin{gather}
\label{eqMeixnerR}
\beta   M_n(x;\beta,c) = (\beta+n)M_n(x;\beta+1,c) - n M_{n-1}(x;\beta+1,c).
\end{gather}
The second order linear operator to be considered is $T_M$,
\begin{gather}
\label{eqTMeixner}
(T_M f)(x) :=  \frac{x!}{c^x(\beta)_x} \nabla\left(\frac{(\beta+2)_xc^x}{x!} \Delta f\right)(x).
\end{gather}

We consider the inner product spaces endowed with the inner product
\begin{gather}
\label{eqinnerpMeixner}
\langle f,g\rangle_\beta = \sum_{x=0}^\infty \frac{c^x(\beta)_x}{x!} f(x) \overline{g(x)}.
\end{gather}
The operator $T_M$  is formally self-adjoint with respect to $\langle \cdot,\cdot\rangle_\beta$.

Using \eqref{eqMeixnerR} and \eqref{eqDEMeixner} we see that
\begin{gather*}
T_M M_n(x;\beta,c) = \frac{(\beta+x)(c-1)}{\beta^2c(\beta+1)} \left[n(\beta+n) M_n(x;\beta+1,c)
- n(n-1)  M_{n-1}(x;\beta+1,c)   \right].
\end{gather*}
 Therefore
 \begin{gather*}
\beta^2  \frac{c(\beta+1)}{(c-1)}  \langle M_m(x;\beta,c), T_M M_n(x;\beta,c)\rangle_\beta \\
\quad = \langle  (\beta+m)M_m(x;\beta+1,c) - m M_{m-1}(x;\beta+1,c),  \\
\quad \qquad\qquad n(\beta+n) M_n(x;\beta+1,c)
- n(n-1)  M_{n-1}(x;\beta+1,c)\rangle_{\beta+1}\\
\quad= m(\beta+m)^2 h_m(\beta+1)\delta_{m,n} - m(m+1)(\beta+m) h_m(\beta+1)\delta_{m,n-1} \\
\qquad\qquad
\quad{}- m(m-1)(\beta+m-1)h_{m-1}(\beta+1) \delta_{m, n+1}+ m^2(m-1) h_{m-1}(\beta+1) \delta_{m,n},
\end{gather*}
where
\begin{gather*}
h_n(\beta) =  \langle M_n(x;\beta,c),  M_n(x;\beta,c)\rangle_\beta =  \frac{c^{-n}n!}{(\beta)_n(1-c)^\beta}.
\end{gather*}
Since $T_M$ annihilates constants we mod out by the space of constants and let the matrix elements
of $T_M$  be $\{B_{m,n}:m,n \ge 0\}$. Thus
\[
B_{m,n}(M) = \frac{\langle M_{m+1}(x;\beta,c), T_M M_{n+1}(x;\beta,c)\rangle_\beta}
{\sqrt{h_{m+1}(\beta)h_{n+1}(\beta)}}.
\]
In other words
\begin{gather}
c\beta(\beta+1) B_{m,n} = - [ (m+1)(m+\beta+1) + m(m+1)c]   \delta_{m,n}\label{eqMatrMei}\\
\hphantom{c\beta(\beta+1) B_{m,n} =}{}
+ m\sqrt{c(m+1)(m+\beta)}  \delta_{m,n+1} + (m+1)\sqrt{c(m+2)(\beta+m+1)}  \delta_{m,n-1}.\nonumber
\end{gather}
Now scale the energy parameter $E$ by $E = -x/(c\beta(\beta+1))$.  This translates into the monic  three
term recurrence relation
\begin{gather}
xP_m(x) = P_{m+1}(x) +  [ (m+1)(m+\beta+1) + m(m+1)c]P_m(x)\nonumber\\
\hphantom{xP_m(x) =}{}
+ cm^2(m+1)(\beta+m)P_{m-1}(x).\label{eqnewpoli1}
\end{gather}
The polynomials generated by  \eqref{eqnewpoli1} seem to be new.   They give a one parameter
generalization of the continuous dual Hahn polynomials which is dif\/ferent from the Wilson
polynomials. Finding the orthogonality measure of these polynomials remains a challenge. It clear
that the measure is unique and is supported on an unbounded  subset of $[0, \infty)$.  They are
birth and death process polynomials corresponding to birth rates $b_n = (n+1)(\beta+n+1)$ and death
rates $d_n = cn(n+1)$, see~\eqref{eqBandDcondition}.
By~\eqref{eqiii} the corresponding moment problem is determinate. So in this case we do not have
a precise analogue of Proposition~\ref{prop:Laguerecase}, except that the unbounded operator~$T_M$ def\/ined
on the polynomials has a unique self-adjoint extension.

\section{Operators with additional potential}\label{section5}

We consider the case of second order operators which arise from classical orthogonal polynomials.
Let $p_n(x)$ be a  monic family of classical orthogonal polynomials and $T$ a second order operator
such that
\begin{gather}
\label{eqTspectrum}
Tp_n(x) = \lambda_n p_n(x).
\end{gather}
Also assume that the three term recurrence relation for the $p_n$'s is
\begin{gather}
\label{eq3trrstart}
xp_n(x) = p_{n+1}(x) + \alpha_n p_n(x) + \beta_n p_{n-1}(x).
\end{gather}
We now consider the spectral problem
\begin{gather}
\label{eqgenspectral}
(T + \gamma x)\psi(x, E) = E\psi(x,E).
\end{gather}
One can think of $Ty = Ey$ as a free particle problem then \eqref{eqgenspectral} will be a
Schr\"{o}dinger problem with potential $\gamma x$. Let $\mu$ be the orthogonality measure of
$\{p_n(x)\}$ and assume we deal with the case when the polynomials are
dense in $L_2({\mathbb{R}}, \mu)$, which is true with very few exceptions. The orthonormal polynomials are
$\{p_n(x)/\sqrt{\beta_1\beta_2\cdots \beta_n}\,\}$ and form a basis for $L_2({\mathbb{R}}, \mu)$. The matrix element
of $T+\gamma x$ with respect to this basis are
\begin{gather*}
B_{m,n} = \frac{1}{\prod\limits_{j=1}^m\sqrt{\beta_j}}  \frac{1}{\prod\limits_{k=1}^n\sqrt{\beta_k}}
 \int_{\mathbb{R}} p_m(x) (T+\gamma x) p_n(x) d\mu(x).
 \end{gather*}
Clearly
\[
B_{m,n} = (\lambda_n+ \gamma \alpha_n) \delta_{m,n}+ \gamma \sqrt{\beta_m}\delta_{m,n+1} + \gamma  \sqrt{\beta_{m+1}}\delta_{m,n-1}.
\]
The matrix $B= \{B_{m,n}\}$ is tridiagonal and generates the  monic orthogonal polynomials via
$\{\phi_n(E)\}$
\begin{gather}
\phi_0(E) =1, \qquad  \phi_1(E) = E - \lambda_0 - \gamma \alpha_0, \nonumber\\
\phi_{n+1}(E) = (E -\lambda_n -\gamma\alpha_n) \phi_n(E) - \gamma^2\beta_n\phi_{n-1}(E).\label{eqdefphi}
\end{gather}
We can still scale $E$ by $E = \xi(x - \eta)$ and introduce additional parameters to help identify
the polynomials as known ones. Thus we let $\psi_n(x) = \xi^{-n}\phi_n(E)$ and transform~\eqref{eqdefphi} to
\begin{gather}
\psi_0(x) =1, \qquad  \psi_1(x) = x -\eta - (\lambda_0 + \gamma \alpha_0)/\xi, \nonumber\\
\psi_{n+1}(x) = [x -\eta -(\lambda_n +\gamma\alpha_n)/\xi] \psi_n(x) - (\gamma/\xi)^2\beta_n\psi_{n-1}(x).\label{eqdefpsi}
\end{gather}
The importance of this scaling  will be made clear in the examples.

Recall that \eqref{eq3trrstart} generates a birth and death process polynomials if there are sequences
$\{b_n\}$ and $\{d_n\}$ such that
\begin{gather}
\label{eqBandDcondition}
\alpha_n = b_n+d_n, \qquad \beta_n = d_n b_{n-1},
\end{gather}
and for $n > 0$, $b_{n-1} >0$ and $d_n >0$, with $d_0 \ge 0$. One can represent the transition
probability of  going  from state $m$ to state $n$ in time $t$ as the Laplace transform of the
product of two orthogonal polynomials and their orthogonality measure. For details, and additional
information and references see
\cite[Chapter~5]{Ism} and the survey article~\cite{Ism:Let:Mas:Val}.

The $b_n$'s and $d_n$'s are birth and death rates at state  (population) $n$. In the case of birth
and death processes with absorption (killing)  Karlin and Tavar\'{e}~\cite{Kar:Tav} showed that
the corresponding orthogonal polynomials satisfy  \eqref{eq3trrstart} where
\eqref{eqBandDcondition} is now replaced by
\begin{gather*}
\alpha_n = b_n+d_n+c_n, \qquad \beta_n = d_n b_{n-1},
\end{gather*}
where $c_n$ is the absorption rate at state~$n$. This leads  to the following remark.

\begin{Remark}
\label{b/d/a}
Assume that $T$ is a positive  linear operator   and~\eqref{eqTspectrum} holds where $\{p_n\}$ are
birth and death process polynomials with birth and death rates $\{b_n\}$ and $\{d_n\}$,
respectively. Then the orthogonal polynomials in~\eqref{eqdefpsi} with $\xi=\gamma$ which arise  in the
tridiagonalization of $T+\gamma x$ are polynomials associated with a birth and death process  with
absorption where the birth and death rates~$\{b_n\}$ and~$\{d_n\}$, respectively and the absorption
rates are $\{\lambda_n/\gamma\}$.
\end{Remark}

The phenomena described in Remark~\ref{b/d/a} seem to be related to shape invariance and related
topics in  Discrete Quantum Mechanics recently developed by R.~Sasaki and his coauthors, see the
recent survey~\cite{Oka:Sas}.

\begin{Example}[Laguerre polynomials]\label{example1}
 In this case $T$ is as on the left-hand side
of \eqref{eqLag2ndorder} and
\begin{gather}
\lambda_n = -n, \qquad \alpha_n = 2n+\alpha+1,\qquad  \beta_n = n(n+\alpha), \qquad
p_n(x) = (-1)^nn!L_n^{(\alpha)}(x).
\notag
\end{gather}
The recursion in \eqref{eqdefpsi} is
\begin{gather}
\label{eqpsiLaguerre}
\psi_{n+1}(x) = \left[x-\eta  + \frac{n(1-2\gamma)}{\xi}  - \frac{\gamma(\alpha+1)}{\xi}\right] \psi_n(x) - n(n+\alpha)\frac{\gamma^2}{\xi^2}\psi_{n-1}(x),
\end{gather}
with  $E = \xi(x-\eta)$.  When $\gamma = 1/4$ we take $\xi = 1/4, \eta = -2\alpha-2$. This identif\/ies
$\psi_n(x)$ as $(-1)^nL_n^{(\alpha)}(x)$. Hence the spectrum is purely continuous and is given by $x \ge0$, that is $E \ge (\alpha+1)/2$.  The absolutely continuous component is given by
\begin{gather*}
\mu^\prime(E) =  4^{\alpha+1}\frac{\exp\left(2(\alpha+1)\right)}{\Gamma(\alpha+1)} [E - (\alpha+1)/2]^\alpha e^{-4E},
\qquad E \in [(\alpha+1)/2, \infty).
\end{gather*}

We next  assume $\gamma >1/4$ and compare \eqref{eqpsiLaguerre} with
 the following monic form of \eqref{eqMeixPoll3trr}
\[
\psi_{n+1}(x) = [x-  (n+\lambda)\cot\phi]  \psi_n(x)- \frac{n(n+2\lambda-1)}{4\sin^2\phi} \psi_{n-1}(x).
\]
We make the parameter identif\/ication
\begin{gather}
\gamma = \frac{1}{4}\sec^2(\phi/2), \qquad \xi = -\tan(\phi/2), \qquad \lambda = (\alpha+1)/2,
\qquad  \eta= \frac{\alpha+1}{2}\cot (\phi/2).
\label{eqparclarge}
\end{gather}
With this choice of parameters we identify the $\psi$'s as Meixner--Pollaczek polynomials. Indeed
$\psi_n(x) = P_n^{(\lambda)}(x;\phi)$  where
\begin{gather}
\gamma = (1+\xi^2)/4, \qquad \lambda = (\alpha+1)/2, \qquad  \eta = - \frac{\alpha+1}{2\xi},
\end{gather}
and
\begin{gather}
\phi_n(E) = \xi^{-n}P_n^{(\lambda)}(\eta + E/\xi).
\end{gather}
The spectral  measure $\mu$ is absolutely continuous and  when normalized to have a total mass
$1$, its Radon--Nikodym derivative is
\begin{gather*}
 \left(\frac{2\xi}{1+\xi^2}\right)^{\alpha+1} \frac{\exp(2\phi-\pi)x)}{\pi \xi  \Gamma(\alpha+1)}
\left|\Gamma(ix+(\alpha+1)/2)\right|^2.
\end{gather*}
We now consider the case $0 <  \gamma < 1/4$. We identify \eqref{eqpsiLaguerre} with the monic form
of \eqref{eqMei3trr}, namely
\[
y_{n+1}(x) = \left[x- \frac{n - c(\beta+n)}{1-c}\right]   y_n(x)
- \frac{cn(n+\beta-1)}{(1-c)^2}y_{n-1}(x).
\]
This is done through the parameter identif\/ication
\begin{gather*}
\gamma = \frac{\sqrt{c}}{(1+ \sqrt{c})^2}, \qquad \xi = - \frac{1-\sqrt{c}}{1+\sqrt{c}}, \qquad
\beta = {\alpha+1}, \qquad \eta = \frac{\beta \sqrt{c}}{1-\sqrt{c}}.
\end{gather*}
It is clear from \eqref{eqMei3trr} that the $y_n$'s are monic  Meixner polynomials.

Note that such a division also occurs in the spectral decomposition of suitable elements in the Lie
algebra $\mathfrak{su}(1,1)$ in the discrete series representations, see~\cite{KoelVdJ, MassR}. So one can ask for Lie algebraic interpretations along the lines of~\cite{Mill},
see~\cite{Ojha} for a related result.
\end{Example}

\begin{Example}[Ultraspherical polynomials]\label{example2}
In this case \cite[(4.5.8)]{Ism}
\[
T = (1-x^2)^{-\nu+1/2}\frac{d}{dx}\left((1-x^2)^{\nu+1/2}\frac{d}{dx}\right), \qquad
\lambda_n =- n(n+2\nu).
\]
The coef\/f\/icients in the monic form of the three term recurrence relation are,
see \cite[(1.8.18)]{Koe:Swa},
\[
\alpha_n =0,\qquad  \beta_n = \frac{n(n+2\nu-1)}{4(n+\nu)(n+ \nu-1)}.
\]
Thus the recursion in \eqref{eqdefpsi} becomes
\begin{gather}
\label{eqAl-Mi}
\psi_{n+1}(x) = \big[x -\eta + \xi^{-1} n(n+2\nu)\big]\psi_n(x) - \frac{\gamma^2   \xi^{-2}n(n+2\nu-1)}{4(n+\nu)(n+\nu-1)}\psi_{n-1}(x).
\end{gather}
We do not know the orthogonality measure of the polynomials in~\eqref{eqAl-Mi}. The special case
$\nu =1$ appeared earlier in the work of Alhaidari and Bahlouli~\cite[(3.8)]{Alh:Bah} where they
applied the $J$-matrix method to quantum model whose potential is an inf\/inite potential well with
sinusoidal bottom.  The same case also appeared in the work~\cite{Goh:Mic} by Goh and Micchelli  on
certain aspects of the uncertainty principle.   Determining the orthogonality measure of these
polynomials will be very useful.

The parameter $\eta$ can be absorbed in $x$, hence we assume $\eta =0$. In the notation of
\eqref{eqorthnormal}
\[
b_n = -\frac{n(n+2\nu)}{\xi}, \qquad a_n^2 =  \frac{\gamma^2 n(n+2\nu-1)}{4 \xi^2(n+\nu)(n+\nu-1)}.
\]
In the case $\xi>0$,  since $ b_n < 0$ and $b_n - b_{n-1} <0$,  Theorem~\ref{IsmLi} implies that the
smallest zero of~$p_n(x)$ is approximately $\frac{1}{2} (b_n+b_{n-1}) - \frac{1}{2}|b_n-b_{n-1}|$.
Hence
\begin{gather*}
x_{n,n} =  -\xi^{-1}n(n+2\nu) + {\mathcal{O}}(1).
\end{gather*}
On the other hand $a_n$ is monotone decreasing if $\nu \ge 1$ or $-1/2 <\nu \le 0$ and monotone increasing if $0 < \nu \le 1$.  Using
\begin{gather*}
\frac{1}{2} (b_n+b_{n-1}) \pm \frac{1}{2} \sqrt{(b_n-b_{n-1})^2 + 16a_n^2} \\
\qquad{} \leq
 \frac{1}{2} (b_n+b_{n-1}) + \frac{1}{2} |b_n-b_{n-1}| + 2 a_n \le 2 \max \{a_1, a_\infty\},
\end{gather*}
where $a_\infty= \lim\limits_{n\to\infty} a_n$.
Therefore
\begin{gather*}
x_{n,1} < 2\max \{a_1, a_\infty\}.
\end{gather*}
Thus the spectrum is unbounded below and is contained in $(-\infty, 2\max \{a_1, a_\infty\})$.  It
is important to note that $p_1(0)=0$, hence the right end point of the spectrum, being
$\lim\limits_{n\to\infty} x_{n,1}$ is positive. The case $\nu =1$ is the case when our starting point is
the Chebyshev polynomials of the second kind. In this case $a_n = a_1$ for all $n$ and the largest
zero of~$p_2(x)$ is
$-\frac{3}{2} + \sqrt{a_1^2 + \frac{9}{4}}>0$.
\end{Example}

\begin{Example}[$q$-Ultraspherical polynomials]\label{example3}
The weight function is supported on $[-1, 1]$ and is given by \cite[\S~13.2]{Ism}
\begin{gather*}
w(x; \beta) dx
= \frac{(e^{2i\theta}, e^{-2i\theta};q)_\infty}{(\beta e^{2i\theta}, \beta e^{-2i\theta};q)_\infty}  d\theta,
\qquad x = \cos \theta, \qquad  \beta < 1.
\end{gather*}
The second order operator is  \cite[\S~13.2]{Ism}
\begin{gather*}
T = \frac{1}{w(x; \beta)} \mathcal{D}_q \left[w(x; q\beta) \mathcal{D}_q \right], \qquad \lambda_n =- \frac{4q^{1-n}}{(1-q)^2}(1-q^n)\big(1-\beta^2q^n\big).
\end{gather*}
In this case
\begin{gather*}
\alpha_n =0, \qquad \beta_n = \frac{(1-q^n)(1-\beta^2 q^{n-1})}{4(1-\beta q^n)(1-\beta q^{n-1})},
\end{gather*}
and   the recurrence relation  in~\eqref{eqdefpsi} gives
\begin{gather}
\psi_{n+1}(x) = \left[x-\eta  + \frac{4q^{1-n}\xi^{-1}}{(1-q)^2}(1-q^n)(1-\beta^2q^n)\right] \psi_n(x)\nonumber\\
 \hphantom{\psi_{n+1}(x) =}{}
 -  \frac{\gamma^2  (1-q^n)(1-\beta^2 q^{n-1})}{4 \xi^2(1-\beta q^n)(1-\beta q^{n-1})} \psi_{n-1}(x).
\label{eqqultra}
\end{gather}

 It is clear $\eta$ can be absorbed in~$x$ so we may assume $\eta =0$. We do not know any explicit
 formulas for the above polynomials. It is clear that they are orthogonal on an unbounded set and that
 Condition~\eqref{eqi} is satisf\/ied,  hence the orthogonality measure is unique.  As in Example~\ref{example2}
 we can show the the spectrum is bounded above and unbounded below and estimate the largest and
 smallest zeros of $p_n(x)$. In the present case
 \begin{gather*}
 b_n = -  \frac{4q^{1-n}\xi^{-1}}{(1-q)^2}(1-q^n)\big(1-\beta^2q^n\big), \qquad a_n^2 =  \frac{\gamma^2  (1-q^n)(1-\beta^2 q^{n-1})}{4 \xi^2(1-\beta q^n)(1-\beta q^{n-1})}.
 \end{gather*}
 Here again the $b_n$'s are negative and decreasing in $n$ for $\xi>0$.  A simple calculation shows that~$a_n$
 increases with~$n$ if $0 < \beta < q$ and decreases with $n$ if $q < \beta < 1$. Thus
 \[
 A := \max \{a_n: n =1, 2, \dots\} = \begin{cases} a_\infty &   \textup{if} \  0 < \beta <q, \\
 a_1  &   \textup{if} \  q < \beta <1.
 \end{cases}
 \]
 Therefore Theorem~\ref{IsmLi}  shows that the smallest zero $x_{n,n}$ satisf\/ies
 \begin{gather*}
 x_{n,n} > \frac{1}{2}(b_n+b_{n-1}) + \frac{1}{2}(b_n - b_{n-1}) -2a_n = b_n -2a_n> b_n -2A.
 \end{gather*}
 Indeed $x_{n,n} = b_n +{\mathcal{O}}(1)$.  To determine the other end of the spectrum note that
 \[
 \sqrt{(b_n-b_{n-1})^2 + 16 a_n^2}  <  |b_n-b_{n-1}| + 4 a_n = b_{n-1} -b_n   + 4 a_n.
 \]
Thus
\begin{gather*}
\frac{1}{2}(b_n+ b_{n-1}) + \frac{1}{2} \sqrt{(b_n-b_{n-1})^2 + 16 a_n^2}    \\
\qquad {}
< \frac{1}{2}(b_n+ b_{n-1}) +
\frac{1}{2}( b_{n-1}-b_n) + 2 a_n = b_{n-1} + 2 a_n \le 2A, \qquad n >0.
\end{gather*}
 Consequently the largest zeros $x_{n,1}$ is $<2A$. Therefore the spectrum of $T+\gamma x$ is
 unbounded below and is contained in $(-\infty, 2A]$.
\end{Example}

\begin{Example}[Chebyshev polynomials]\label{example4}
The Chebyshev polynomials of the f\/irst and second kinds are special ultraspherical polynomials and
special
$q$-ultraspherical polynomials as well.  We will only discuss the polynomials $\{U_n(x)\}$ but the
reader can easily write down the corresponding  formulas for the polynomials~$\{T_n(x)\}$.

The $U_n$'s correspond to $\nu =1$ of \eqref{eqAl-Mi} and the case $\beta = q$ of~\eqref{eqqultra}.
Thus we are led to  the following systems of orthogonal polynomials
\begin{gather*}
r_{n+1}(x) = \big[x -\eta + \xi^{-1} n(n+2)\big]r_n(x) - \frac{\gamma^2}{4   \xi^{2}}r_{n-1}(x), \\
s_{n+1}(x) = \left[x-\eta  + \frac{4q^{1-n}\xi^{-1}}{(1-q)^2}(1-q^n)\big(1- q^{n+2}\big)\right] s_n(x)
 -  \frac{\gamma^2  }{4 \xi^2} s_{n-1}(x).
\end{gather*}
Here again we do not know any explicit representations or orthogonality measures for the
polynomials $\{r_n(x)\}$ or $\{s_n(x)\}$.  Again condition~\eqref{eqi} is satisf\/ied for
$\{r_n(x)\}$ and $\{s_n(x)\}$. Therefore the orthogonality  measures of both   families of
polynomials   unique.
\end{Example}

\section{Adding a linear  potential}\label{section6}
In this section we  yet have a variation on the problems of potential introduced at the
beginning  of  Section~\ref{section5}.  We start with \eqref{eqTspectrum} where the
eigenfunctions  satisfy~\eqref{eq3trrstart}. We then consider
the Schr\"{o}dinger operator
\begin{gather} \label{eq:61}
S = (x+c)T+ \gamma x.
\end{gather}
We  illustrate this idea by considering the operators $T_L$ and $T_M$ for the Laguerre and Meixner
polynomials def\/ined in  \eqref{eq3.1} and \eqref{eqTMeixner}, respectively.

{\bf The Laguerre case.}
 Here we take $c =0$ and
\[
T = \frac{e^x}{ x^{\alpha+1}} \frac{d}{dx} \left[x^{\alpha+2}e^{-x} \frac{d}{dx}\right].
 \]
With the notation in \eqref{eq3.1} we have
\begin{gather*}
S = T_L + \gamma x.
\end{gather*}
We use the inner product \eqref{eqinnerprodL}  and our weighted
$L_2$ space is $L_2(0, \infty, x^\alpha e^{-x})$. The matrix elements
$S_{m,n}$ are
\begin{gather*}
\sqrt{\frac{m!n!}{\Gamma(\alpha + m +1) \Gamma(\alpha +n +1)}}\\
 \qquad{}\times
\left[ \int_0^\infty  L_m^{(\alpha)}(x)
\frac{d}{dx} \left\{x^{\alpha+2}e^{-x} \frac{d}{dx} L_n^{(\alpha)}(x)\right\} dx
+ \gamma\int_0^\infty L_m^{(\alpha)}(x) xL_n^{(\alpha)}(x) x^\alpha e^{-x} dx\right].
\end{gather*}
Using the recurrence relation
\begin{gather*}
xL_n^{(\alpha)}(x) = -(n+1) L_{n+1}^{(\alpha)}(x) - (n+\alpha) L_{n-1}^{(\alpha)}(x)
+ (2n+\alpha+1)L_n^{(\alpha)}(x),
\end{gather*}
and the calculation of the matrix elements in   \eqref{eq*} we f\/ind that
\begin{gather*}
S_{m,n} = [\gamma(2m+\alpha+1)-m(\alpha+2m)] \delta_{m,n}   \\
\hphantom{S_{m,n} =}{}  +(m-\gamma) \sqrt{(m+1)(m+\alpha+1)}   \delta_{m,n-1} +
(m-1-\gamma)\sqrt{m(m+\alpha)}   \delta_{m,n+1}.
\end{gather*}
The null space of $S$ is trivial so there is no need to mod out by the null space as we did in
Sections~\ref{section3} and~\ref{section4}.  The monic polynomials $\{p_n(E)\}$  which arise through tridiagonalization are
generated by
 $p_0(E)=1$,   $p_1(E)= E- \gamma(\alpha+1)$, and
\begin{gather}
Ep_n(E) = p_{n+1}(E) + n(n+\alpha) (n-1-\gamma)^2p_{n-1}(E)\nonumber\\
\hphantom{Ep_n(E) =}{}
+ [\gamma(2n+\alpha+1)-n(\alpha+2n)]p_n(E).\label{eq5.25}
\end{gather}
The polynomials in \eqref{eq5.25} form a two parameter subfamily of the continuous dual Hahn
polynomials \cite[\S~1.3]{Koe:Swa}  with the parameters
\[
a = -\gamma -(\alpha+1)/2, \qquad b = c = (\alpha+1)/2, \qquad E = - x -\frac{1}{4}(\alpha+1)^2.
\]
In the above analysis we  assumed $\alpha > -1$, hence $b = c >0$. If $\gamma < - (\alpha+1)/2$ then $S$
has purely a continuous spectrum supported on $(-\infty, -(\alpha+1)^2/4]$, see~\eqref{eq:orthorelcontHahn}.
If $\gamma > -(\alpha+1)/2$ there is also discrete spectrum, and we obtain as in Section~\ref{section3}
the following proposition.

\begin{prop}\label{prop:Laguerrewithlinearpotential}\sloppy
The unbounded operator $S$ acting on the subspace of polynomials in
  $L_2(0, \infty, x^\alpha e^{-x})$ is essentially self-adjoint. The spectrum of the closure has
an absolutely continuous part $(-\infty, -\frac14(\alpha+1)^2]$. The discrete spectrum consists
of $\{E_k\mid k\in \{0,\dots, M\}\}$, where $M = \max \{k\in{\mathbb{N}}\mid  k -\gamma- (1+\alpha)/2 < 0\}$
and $E_k = (k-\gamma)(k-\gamma-\alpha-1)\in (-\frac14(\alpha+1)^2, \gamma (\gamma+\alpha+1)]$.
\end{prop}

\begin{Remark} The case $\gamma=0$ of Proposition \ref{prop:Laguerrewithlinearpotential} reduces to
Proposition~\ref{prop:Laguerecase} using a similar reduction as in \cite[Remark~3.4]{Ism:Koe2}. Note that
in this case $E_0=0$.
\end{Remark}

 \begin{Remark}
 In this case $S= T_L + \gamma x$ can be written as $S = x \left(D^{\alpha+1} + \gamma\right)$, where
 $D^{\alpha+1}$ is the second order dif\/ferential operator
 \[
 D^{\alpha+1} = \frac{1}{x^{\alpha+1}e^{-x}} \frac{d}{dx} \left[x^{\alpha+2}e^{-x} \frac{d}{dx}\right],
 \]
 for which $D^{\alpha+1}L_n^{(\alpha+1)}(x) = -n L_n^{(\alpha+1)}(x)$.  This shows that $S$ is of the type considered in
 \cite{Ism:Koe2}.
 \end{Remark}

  {\bf The Meixner case.}
With the notation in \eqref{eqTMeixner} we let
\begin{gather*}
\tilde{S}   = T_M + \frac{(1-c) \gamma}{c\beta (\beta+1)}  x.
\end{gather*}
This corresponds to $c=0$ in \eqref{eq:61}, and $T_M$ as in \eqref{eqTMeixner}.
We use the inner product \eqref{eqinnerpMeixner}  and our space is now $L_2$  weighted with
the orthogonality measure of the Meixner polynomials with parameters $c$ and $\beta$.
 The matrix elements
$\tilde{S}_{m,n}$ are
\begin{gather*}
\sqrt{\frac{(\beta)_m(\beta)_n(1-c)^{2\beta}}{c^{-m-n}m!n!}}
\left[ \,\sum_{x=0}^\infty  M_m(x; \beta,c)
\nabla\left(\frac{(\beta+2)_xc^x}{x!} \Delta M_n(x; \beta,c) \right) \right.\\
\left. \qquad{}
+ \frac{(1-c)\gamma}{c\beta (\beta+1)}   \sum_{x=0}^\infty \frac{c^x (\beta)_x}{x!}  M_m(x; \beta,c)  xM_n(x; \beta,c) \right].
\end{gather*}
We already calculated the matrix elements of $T_M$ in~\eqref{eqMatrMei}, but we must replace $m$, $n$
by $m-1$, $n-1$, respectively. Using the recurrence relation~\eqref{eqMei3trr}  we then compute the
matrix elements of a constant times~$x$. This leads to
\begin{gather*}
c\beta(\beta+1) \tilde{S}_{m,n} = - [ m(m+\beta) + m(m-1)c +\gamma m+c\gamma(\beta+m)]   \delta_{m,n} \\
\qquad{}+ (m-1-\gamma)\sqrt{c m(m+\beta-1)}\, \delta_{m,n+1} + (m-\gamma) \sqrt{c(m+1)(\beta+m)}  \delta_{m,n-1}.
\end{gather*}
Therefore the corresponding orthonormal polynomials $p_n(E)$ satisfy
the three term recurrence relation
\begin{gather}
c\beta(\beta+1)Ep_m(E) = - [ m(m+\beta) + m(m-1)c +\gamma m+c\gamma(\beta+m)] p_m(E) \nonumber\\
\!\qquad{}+ (m-\gamma) \sqrt{c(m+1)(\beta+m)}p_{m+1}(E) + (m-1-\gamma)\sqrt{c m(m+\beta-1)}p_{m-1}(E).\!\!\!\label{eqNewMexpol}
\end{gather}
We assume $0<c<1$, $\beta>0$, since we are dealing with the Meixner polynomials. In order to have
\eqref{eqNewMexpol} satisfy the conditions for an orthonormal polynomial system we assume
\mbox{$\gamma<0$}.
 The polynomials generated by~\eqref{eqNewMexpol} seem to be new. Here again we have neither
 explicit representations or generating functions, nor do we know their orthogonality measure. One
 can say however that their orthogonality measure is unique since condition~\eqref{eqii} is clearly
 satisf\/ied for suf\/f\/iciently large~$n$.  Using Theorem~\ref{IsmLi} and some estimates we see that
 the support of the orthogonality measure is contained in $(-\infty, a]$ for some $a>0$.

\section[The Al-Salam-Chihara polynomials]{The Al-Salam--Chihara polynomials}\label{section7}

Recall that $w$ is def\/ined in \eqref{eqweightAC}.  The generating function \eqref{eqgfAC} implies
\begin{gather}
\label{eqContiAC}
(1-t_1t_2)p_n(x; t_1,t_2)
= (1-t_1t_2q^n) p_n(x; t_1,qt_2)  -
t_1 t_2 (1-q^n) p_{n-1}(x; t_1, qt_2).
\end{gather}

We will f\/irst consider the case when  the  operator  $T$ is
\begin{gather*}
L := \frac{1}{w(x; t_1,t_2)} \mathcal{D}_q \left[w\big(x ; q^{1/2}  t_1,q^{3/2} t_2\big)  \mathcal{D}_q \right].
\end{gather*}
Apply (15.1.6) and (12.2.2) in~\cite{Ism} to see that
\begin{gather}
\label{eqAClowering}
\mathcal{D}_q p_n(x; t_1, t_2) = \frac{2t_1 q^{1-n}(1-q^n)}{(1-q)(1-t_1t_2)}
p_{n-1}\big(x; t_1\sqrt{q}, t_2\sqrt{q}\big).
\end{gather}
Moreover
\begin{gather}
\label{eqACraise}
  \mathcal{D}_q \left[ w\big(x, \sqrt{q}t_1, \sqrt{q}t_2\big) p_{n-1}\big(x; \sqrt{q}t_1, \sqrt{q}t_2\big) \right]
= \frac{2(1-t_1t_2)}{t_1(q-1)}  w(x, t_1, t_2) p_n(x;  t_1, t_2).
\end{gather}
 For real $t_1$, $t_2$ with $|t_1|, |t_2|<1$ the Al-Salam--Chihara polynomials  $\{p_n(x; t_1, t_2)\}$ satisfy the orthogonality relation
 \cite[(15.1.5)]{Ism},
 \begin{gather}
 \label{eqorthAC}
h_n(t_1,t_2) \delta_{m,n} =
 \int_{-1}^1   p_m(x; t_1,t_2) p_n(x; t_1,t_2) w(x; t_1,t_2) dx, \\
 h_n(t_1,t_2) = \frac{2\pi (q;q)_nt_1^{2n}}{(q, t_1t_2;q)_\infty(t_1t_2;q)_n},\nonumber
 \end{gather}
 and  are complete in $L^2(-1,1;w(x; t_1, t_2)dx)$.  In view of \eqref{eqorthAC} the orthonormal Al-Salam--Chihara polynomials are
\begin{gather*}
\tilde{p}_n(x;t_1, t_2) =  \sqrt{\frac{(q,t_1t_2;q)_\infty(t_1t_2;q)_n}{2\pi t_1^{2n}(q;q)_n}}   p_n(x; t_1, t_2).
\end{gather*}

\begin{lem}
Let $\{A_{m,n}(AC)\}$ be the  matrix elements of $L$  in the basis $\{\tilde{p}_n(x;t_1, t_2)\}$. Then
\begin{gather*}
A_{m,n} = -\frac{4q^{1-m}(1-q^m)}{(1-q)^2}[1-t_1t_2q^m + t_2^2(q-q^m)]  \delta_{m,n}\\
\hphantom{A_{m,n} =}{} + \frac{4t_2q^{1-m}(1-q^m)}{(1-q)^2}\sqrt{(1-q^{m+1})(1-t_1t_2q^m)} \delta_{m,n-1} \\
\hphantom{A_{m,n} =}{}+ \frac{4t_2q^{2-m}(1-q^{m-1})}{(1-q)^2}  \sqrt{(1-q^{m})(1-t_1t_2q^{m-1})}
 \delta_{m,n+1}.
\end{gather*}
\end{lem}

\begin{proof} Clearly equation \eqref{eqAClowering} implies
\begin{gather*}
\int_{-1}^1 p_m(x; t_1,t_2)w(x; t_1, t_2)Lp_n(x;t_1,t_2) dx \\
\qquad{}= \frac{2t_1(1-q^n)q^{1-n}}{(1-q)(1-t_1t_2)}   \int_{-1}^1 p_m(x; t_1, t_2)
\mathcal{D}_q \left[w\big(x ; q^{1/2}t_1,q^{3/2} t_2\big) p_{n-1}(x;  t_1 \sqrt{q}, t_2\sqrt{q}) \right] dx.
\end{gather*}
In view of \eqref{eqContiAC}  the integrand in the last step  is
\begin{gather*}
\left[\frac{(1-t_1t_2q^m)}{1-t_1t_2}  p_m(x; t_1,qt_2)  -
\frac{t_1 t_2 (1-q^m)}{1-t_1t_2}  p_{m-1}(x; t_1, qt_2) \right] \\
\qquad{}\times \mathcal{D}_q \left[w\big(x ; q^{1/2}t_1,q^{3/2} t_2\big)\left\{\frac{(1-t_1t_2q^{n})}{1-q t_1t_2}  p_{n-1}\big(x; \sqrt{q} t_1,q^{3/2}t_2\big)
\right. \right. \\ \left.\left.
 \qquad\quad{} - \frac{t_1 t_2 (q-q^n)}{1-qt_1t_2}  p_{n-2}\big(x; \sqrt{q} t_1, q^{3/2}t_2\big)\right\}\right].
\end{gather*}
Applying  \eqref{eqACraise} we see that the quantity after the $\times$ is
\[
\frac{2}{t_1(q-1)} w(x; t_1, qt_2)
\big[(1-t_1t_2q^n)p_n(x;t_1,qt_2) - t_1t_2(q-q^n)p_{n-1}(x;t_1,qt_2)\big].
\]
Therefore
\begin{gather*}
-\frac{(1-q)^2(1-t_1t_2)^2} {4 (1-q^n)q^{1-n}}\int_{-1}^1 p_m(x; t_1,t_2)w(x; t_1, t_2)Lp_n(x;t_1,t_2) dx
  \\
\qquad{} =  \int_{-1}^1w(x,t_1,qt_2)  [(1-t_1t_2q^m)  p_m(x; t_1,qt_2)  -
t_1 t_2 (1-q^m)  p_{m-1}(x; t_1, qt_2) ]  \\
\qquad\quad{} \times  [(1-t_1t_2q^n)p_n(x;t_1,qt_2) - t_1t_2(q-q^n)p_{n-1}(x;t_1,qt_2)] dx \\
\qquad{}= [(1-t_1t_2q^m)^2h_m(t_1,qt_2)+ t_1^2 t_2^2 (1-q^m) (q-q^m)h_{m-1}(t_1, qt_2)]\delta_{m,n}\\
 \qquad\quad{} - qt_1t_2(1-q^m) (1-t_1t_2q^m) h_m(t_1,qt_2)\delta_{m, n-1} \\
 \qquad\quad{} - t_1t_2(1-q^m)(1-t_1t_2q^{m-1})h_{m-1}(t_1,qt_2) \delta_{m,n+1}.
\end{gather*}
The result follows since
\begin{gather*}
A_{m,n} = \frac{1}{\sqrt{h_m(t_1,t_2)h_n(t_1,t_2)}}
\int_{-1}^1 p_m(x; t_1,t_2)w(x; t_1, t_2)Lp_n(x;t_1,t_2) dx. \tag*{\qed}
\end{gather*}
\renewcommand{\qed}{}
\end{proof}

The monic orthogonal polynomials generated by the matrix $\frac{4}{(1-q)^{2}}A$, with $A= \{A_{m,n}(AC)\}$
satisfy
\begin{gather}
P_{n+1}(x) = [x- q^{-n}(1-q^{n+1})[1-t_1t_2 q^{n+1}+t_2^2q(1-q^n)] P_n(x) \nonumber\\
\hphantom{P_{n+1}(x) =}{} - t_2^2q^{2-2n}(1-q^{n+1})(1-q^n)^2(1-t_1t_2q^{n})P_{n-1}(x).\label{eqnewpol}
\end{gather}

 Note that the recurrence relation \eqref{eqnewpol} is invariant under $q \to 1/q$ after scaling
 and renaming the parameters.  The recurrence coef\/f\/icients grow exponentially, and by \cite[Theorem~VII.1.5]{Ber} we can easily check that the moment problem corresponding problem does not have a
 unique solution (indeterminate).  Hence $L$ with domain the polynomials is not essentially
self-adjoint.
Nothing is known about the explicit formulas of the  polynomials
 generated by~\eqref{eqnewpol} or any of their orthogonality measures.

{\bf The $\boldsymbol{q^{-1}}$-Hermite  polynomials.}
We now study the $q^{-1}$-Hermite polynomials of Askey \cite{Ask}, Ismail and Masson \cite{Ism:Mas}.   They are generated by $h_0(x|q), h_1(x|q) = 2x$, and
\begin{gather*}
 h_{n+1}(x|q) = 2xh_n(x|q) - q^{-n}(1-q^n)h_{n-1}(x|q).
\end{gather*}
Here we use the parametrization $x = \sinh \xi$.   Recall the def\/initions \cite[Chapter~21]{Ism}, \cite{Ism93}
\begin{gather*}
f(x) = \breve{f}\big(e^\xi\big),  \qquad (\mathcal{D}_q f) := \frac{\breve{f}(q^{1/2}e^\xi)- \breve{f}(q^{-1/2}e^\xi)}{(q^{1/2}- q^{-1/2})(e^\xi+e^{-\xi})/2}, \\
(\mathcal{A}_q f)(x) = \frac{1}{2} \big[\breve{f}\big(q^{1/2}e^\xi\big)+ \breve{f}\big(q^{-1/2}e^\xi\big)\big].
\end{gather*}
The corresponding moment problem is indeterminate but all the $N$-extremal measures have been determined in
\cite{Ism:Mas}, see also \cite[\S~4]{Chr:Koe} for another proof.
They are purely discrete and are enumerated by a parameter $a \in (q,1)$. The support is $\{x_n(a): n = 0, \pm 1, \pm 2, \dots\}$ and
\begin{gather}
\label{eqNextr}
x_n(a) =\frac{1}{2} [q^{-n}/a -a q^n], \qquad \mu(\{x_n(a)\}) = \frac{a^{4n} q^{n(2n-1)}(1+ a^2 q^{2n})}{(-a^2, -q/a^2,q; q)_\infty},
\end{gather}
where $\mu$ is the corresponding normalized orthogonality measure.  The orthogonality relation is
\begin{gather*}
\int_{\mathbb{R}} h_m(x|q)h_n(x|q) \, d\mu(x) =   q^{-n(n+1)/2} (q;q)_n  \delta_{m,n}.
\end{gather*}

The lowering operator is
\begin{gather*}
\mathcal{D}_q h_{n}(x|q) = \frac{2(1-q^n)}{1-q} q^{(1-n)/2}h_{n-1}(x|q).
\end{gather*}
The second order operator equation satisf\/ied by the $q^{-1}$-Hermite polynomials is \cite{Ism93}
\begin{gather}
\label{eqqinverseOp2}
q^{1/2}(1+2x^2) \mathcal{D}_q^2 y + \frac{4q}{q-1} x \mathcal{A}_q \mathcal{D}_q y = \lambda_n y, \qquad \lambda_n :=  -\frac{4q(1-q^n)}{(1-q)^2}.
\end{gather}
With the measure $\mu$ def\/ined as in \eqref{eqNextr} the matrix elements of the operator
$T_H+\gamma x$ with $T_H$ the operator on the
left  side of \eqref{eqqinverseOp2} on $L_2({\mathbb{R}}, \mu)$ with basis
 $\{\tilde h_n = q^{n(n+1)/4}h_n(x|q)/\sqrt{(q;q)_n} \}$  are given by
 \[
 (T_H+\gamma x) \tilde h_n = \frac{\gamma}{2}q^{-(n+1)/2}\sqrt{1-q^{n+1}}\tilde h_{n+1}
 + \lambda_n \tilde h_n + \frac{\gamma}{2}q^{-n/2}\sqrt{1-q^{n}}\tilde h_{n-1}.
 \]
Therefore the corresponding monic polynomials $\{p_n(E)\}$ are generated by
\begin{gather}
\label{eq7.16}
Ep_n(E) = p_{n+1}(E) -\frac{4q(1-q^n)}{(1-q)^2}p_n(E)
+ \frac{\gamma^2}{4} q^{-n}(1-q^{n}) p_{n-1}(E).
\end{gather}
This is essentially a perturbation of the Jacobi matrix  of the $q^{-1}$-Hermite polynomials by the
diagonal matrix $-\frac{4q(1-q^n)}{(1-q)^2}$.  Apart from the shift $\frac{-4q}{(1-q)^2}$Id, this
is a compact perturbation of the Jacobi matrix for the $q^{-1}$-Hermite polynomials. Using
\cite[Chapter~VII, Theorem~1.5]{Ber} the moment problem corresponding to the orthogonal polynomials
generated by \eqref{eq7.16}
 is indeterminate.

 We now give bounds for the zeros of $p_n(E)$. In the present case
 \[
 a_n^2 = \frac{\gamma^2}{4} q^{-n}(1-q^{n}), \qquad b_n =- \frac{4q(1-q^n)}{(1-q)^2}.
 \]
 Take $\gamma >0$, it is clear that $a_n$ is monotonic increasing while $b_n$ is monotic decreasing.  The use of
 \[
 \sqrt{(b_n-b_{n-1})^2 + 16 a_n^2} < |b_n-b_{n-1}| + 4 a_n,
 \]
 shows that the $x_j$'s and $y_j$'s in \eqref{eqdefxiyi} satisfy
 \begin{gather*}
 x_j <  \frac{1}{2}( b_{j}+ b_{j-1} ) + \frac{1}{2}(b_{j-1}- b_j) + 2 a_j = b_{j-1} + 2a_j < 2a_n, \\
 y_j > \frac{1}{2}( b_{j}+ b_{j-1} ) - \frac{1}{2}(b_{j-1}- b_j) -  2 a_j = b_j -  2 a_j
 \ge -2a_n - \frac{4q}{(1-q)^2}.
 \end{gather*}
 Therefore
 \begin{gather*}
 x_{n,1} <  \gamma q^{-n/2}\sqrt{1-q^{n}},\qquad  x_{n,n} > -  \gamma q^{-n/2}\sqrt{1-q^{n}}  - \frac{4q}{(1-q)^2}.
\end{gather*}

\subsection*{Acknowledgements}
The research of Mourad E.H. Ismail is supported by a Research Grants Council of Hong Kong  under
contract \# 101411 and NPST Program of King Saud
University, Saudi Arabia, 10-MAT 1293-02.
This work was also partially supported by a grant from the `Collaboration Hong Kong~-- Joint
Research Scheme' sponsored by the Netherlands Organisation of Scientif\/ic Research and the Research
Grants Council for Hong Kong (Reference number: 600.649.000.10N007). The work for this paper was done
while both authors visited City University Hong Kong, and we are grateful for the hospitality.

We thank Luc Vinet and Hocine Bahlouli for useful comments and references. We also thank the referees for their very careful reading and for their suggestions and constructive criticism that have improved the paper.

\pdfbookmark[1]{References}{ref}
\LastPageEnding

\end{document}